\documentclass{scrartcl}

\usepackage{ucs}
\usepackage[utf8]{inputenc}

\usepackage{amsmath}
\usepackage{amscd}
\usepackage{amssymb}

\newcommand{\abs}[1]{\left|#1\right|}
\newcommand{\BM}{\mathrm{BM}}
\newcommand{\ch}{\mathrm{ch}}
\newcommand{\colim}{\operatorname*{colim}}
\newcommand{\coloneqq}{\mathrel{:}=}
\newcommand{\gr}{\operatorname{gr}}

\newcommand{\pt}{\mathrm{pt}}
\newcommand{\res}{\operatorname*{res}}
\newcommand{\rk}{\operatorname{rk}}
\newcommand{\sheaf}[1]{\mathcal #1}
\newcommand{\set}[1]{\mathbf #1}
\newcommand{\Set}[1]{\left\{#1\right\}}
\newcommand{\supp}{\operatorname{supp}}

\newcommand{\units}{\times}
\newcommand{\vac}{\left|0\right\rangle}

\title{Equivariant cohomology, symmetric functions and the Hilbert scheme of points on the total space of the invertible sheaf $\sheaf O_{\set P^1}(-2)$ over the projective line}
\author{Marc A.~Nieper-Wißkirchen\thanks{The author likes to thank S.~Boissière for pointing out to him the articles of H.~Nakajima and W.-P.~Li, Zh.~Qin and W.~Wang on Hilbert schemes and Jack symmetric functions and for carefully reading preliminary versions of this manuscript.}}
\date{\today}

\begin{document}

\maketitle

\abstract{
Let $X$ be the quasi-projective symplectic surface that is given by the total space of the invertible sheaf $\sheaf O_{\set P^1}(-2)$ over the projective line. Let $(X^{[n]})_{n \ge 0}$ be the family of Hilbert schemes of points on $X$. 

We give and prove a closed formula expressing any multiplicative characteristic class evaluated on the
$X^{[n]}$ in terms of the standard Fock space description of the cohomology of the $X^{[n]}$. As a side result, we also deduce a formula for the Chern character of the tangent bundles of the $X^{[n]}$.

The results found here are another step towards a complete description of the tangent bundle of the Hilbert scheme of a general quasi-projective surface as the formulas given here yield certain coefficients in that description.

Finally, we also give a closed formula expressing any multiplicative characteristic class evaluated on the tautological bundles $\sheaf O_X^{[n]}$.
}

\tableofcontents

\section{Introduction}

\subsection{Preface}

Let $X$ be a quasi-projective surface over the complex numbers. The Hilbert scheme $X^{[n]}$, $n \in \set N_0$, parametrises zero-dimensional subschemes of $X$ of length $n$. By a result of J.~Fogarty (\cite{smoothness}) it is a smooth (quasi-projective) variety. 

I.~Grojnowski (\cite{hilbg}) and H.~Nakajima (\cite{hilbn}) have shown that the direct sum of the (rational) cohomology spaces of all $X^{[n]}$ carries a natural structure as a Fock space. It is an irreducible representation over a Heisenberg algebra generated by the cohomology ring $H^*(X, \set Q)$ of the surface. Thus every cohomology class in $H^*(X^{[n]}, \set Q)$ has a natural combinatorial description in terms of this Fock space.

In particular, given a multiplicative characteristic class $\phi$ (i.e.~a characteristic class giving rise to a genus, see~\cite{multiplicative}), one can ask for the value of the cohomology class
\[
\phi(X^{[n]}) \coloneqq \phi(T_{X^{[n]}})
\]
in terms of the creation operators of the aforementioned Fock space. It turns out that there is a universal formula, in which the surface $X$ enters only through its characteristic classes, see e.g.~\cite{boissiere} or~\cite{affine}. Precisely, the shape of the universal formula is
\[
\sum_{n \ge 0} \phi(X^{[n]}) = \exp\sum_{\lambda} \left(a^\lambda_\phi \, q_\lambda(1) + b^\lambda_\phi q_\lambda(K_X) + c^\lambda_\phi \, q_\lambda(e_X) + d^\lambda_\phi \, q_\lambda(K^2_X)\right) \vac
\]
for certain coefficients $(a^\lambda_\phi)$, $(b^\lambda_\phi)$, $(c^\lambda_\phi)$ and $(d^\lambda_\phi)$ where $\lambda$ runs through all partitions. To solve the problem, it thus remains to calculate these coefficients in closed form.

This problem, however, is far from being solved. In~\cite{affine}, S.~Boissière and the author have given a closed expression for the $a^{(k)}_\phi$, i.e.~for the first series of coefficients and only for partitions of length one. In this article, we give and prove a formula describing the $a^{(k, l)}_\phi$, i.e.~we have extended the previous results to partitions of length two.

As any characteristic class, also the Chern character can be calculated from the knowledge of all multiplicative classes. Thus, from our results about multiplicative classes, we are able to prove an analoguous result for the Chern character of the tangent bundles of the Hilbert schemes of points on $X$ as a corollary.

Our calculations are enough to fully solve the problem for surfaces $X$ with $H^4(X, \set Q) = 0$ and with vanishing canonical divisor, $K_X = 0$. In this case, all the other coefficients do not contribute. In fact, we are proving our theorem by doing explicit calculations on one such surface, namely the total space of the invertible sheaf $\sheaf O_{\set P^1}(-2)$ on $\set P^1$. On this surface, there is a natural action by the one-dimensional torus with two isolated fixpoints. W.~P.~Li, Zh.~Qin and W.~Wang building on previous work by H.~Nakajima~(\cite{jackn}) and E.~Vasserot~(\cite{vasserot}) showed in~\cite{jack} how to use equivariant cohomology with respect to this action to express Fock space elements of the ordinary cohomology space by equivariant cohomology classes using symmetric functions.

We, in turn, make heavy use of their description, in order to actually perform the calculations outlined there. For the rest of our article, we assume that the reader has access to that article and the
aformentioned article~\cite{affine}.

\subsection{The Theorem.}

\paragraph{The setting.}
Let $X$ be the total space of the line bundle $\sheaf O_{\set P^1}(-2)$ over the complex projective line. We denote the cohomology class of a fibre by $h \in H^2(X, \set Q)$.

For $n \in \set N_0$, let $X^{[n]}$ be the Hilbert scheme of $n$ points on $X$. Note that $X^{[0]} = \Set \pt$. The unit in $H^*(X^{[0]}, \set Q)$ is denoted by $\vac$.

For $\alpha \in H^*(X, \set Q)$, we denote the Fock space creation operators by
\[
q_m(\alpha)\colon H^*(X^{[n]}, \set Q) \to H^{* + \deg \alpha + 2 \, n - 2}(X^{[m + n]}, \set Q)
\]
where $m \in \set N_0$.

\paragraph{Theorem.}
Let $\phi$ be a multiplicative characteristic class over a $\set Q$-algebra $A$, i.e.~there is a power series $f \in 1 + x \, A[[x]]$ with
\[
\phi = \prod_i f(x_i),
\]
where the $x_i$ are the Chern roots. Let $g \in xA[[x]]$ be the compositional inverse of
\[
\frac{x}{f(x) f(-x)}.
\]

Then, the multiplicative class evaluated on the tangent bundle of the $X^{[n]}$ is given by
\[
\sum_{n \ge 0} \phi(X^{[n]}) = \exp\left(\sum_{k \ge 1} a_k \, q_k(1) +
\sum_{\substack{k, l \ge 1}} a_{k, l} \, q_k(h) \, q_l(h)\right) \, \vac,
\]
where the sequences $(a_k)_k$ and $(a_{k, l})_{k, l}$ in $A$ are defined by
\[
\sum_{k \ge 1} k a_k \, x^k = g(x)
\]
and
\[
\sum_{k, l}^\infty a_{k, l} \, x^k y^l = \log\frac{g(x) - g(y)}{(x - y) \, f(g(x) - g(y)) \, f(g(y) - g(x))}.
\]

\paragraph{Corollary.}

The Chern character evaluated on the tangent bundles of the $X^{[n]}$ is given by
\[
\sum_{n \ge 0} \ch(X^{[n]}) = \left(\sum_{k \ge 1} a_k q_k(1) + \sum_{k, l \ge 1} a_{k, l} q_k(h) q_l(h)\right) \,
\exp(q_1(1)) \, \vac,
\]
where the sequences $(a_k)_k$ and $(a_{k, l})_{k, l}$ in $\set Q$ are defined by
\[
\sum_{k \ge 1}^\infty a_k \, x^k = \sum_{m \ge 0} \frac{2}{(2 m + 1)!} \, x^{2 m + 1}.
\]
and
\[
\sum_{k, l \ge 1} a_{k, l} \, x^k \, x^l
= \sum_{m \ge 1} \frac 2 {(2 m)!} \sum_{k + l = 2 m} \left(
1 - (-1)^k \binom{2 m} k
\right) \, x^k y^l.
\]

The proof of the theorem and the proof of the corollary are given in the last section of this article.

\section{Preliminaries}

\subsection{Notation}

\paragraph{Lie groups.}
By $T$ we denote the Lie group $U(1)$, i.e.~the one-dimensional torus.

\paragraph{Coefficients.}

Let $A$ be an integral domain and $f \in A((z))$ a Laurent series. By
\[
[z^n] f(z)
\]
we denote the coefficient of $f(z)$ before $z^n$. This notion shall be extended to multivariate Laurent series. In particular, we have
\[
[z^{-1}] f(z) = \res_z f(z).
\]

\paragraph{Partitions.}
A \emph{partition $\lambda$} is a weakly decreasing sequence $(\lambda_1, \lambda_2, \dots)$ of non-negative integers.
We set
\[
\abs\lambda \coloneqq \sum_i \lambda_i.
\]
If $n = \abs\lambda$, we call $\lambda$ a partition of $n$.
Its \emph{length} is given by
\[
l(\lambda) = \max_i\Set{\lambda_i \neq 0}.
\]
If the length of $\lambda$ is less or equal than $l$, we often write $(\lambda_1, \ldots, \lambda_l)$ instead of 
$(\lambda_1, \ldots, \lambda_l, 0, \ldots)$.

By $D_\lambda$ we denote its Young diagram (see~\cite{symmetric}). For each cell $w$ in $D_\lambda$, let $l(w)$ be the \emph{leg length} of $w$ (number of cells below $w$) and $a(w)$ the \emph{arm length} of $w$ (number of cells right of $w$). Then
\[
h(w) \coloneqq a(w) + l(w) + 1
\]
is called the \emph{hook length} of $w$.

Let $\alpha, \beta \in \set R$. We define the multiset
\[
W_\lambda(\alpha, \beta) \coloneqq \Set{\alpha \, (l(w) + 1) + \beta \, a(w), - \alpha \, l(w) - \beta \, (a(w) + 1) \mid w \in D_\lambda}.
\]
Following~\cite{jack}, we further define
\begin{align*}
c_\lambda(\alpha, \beta) & \coloneqq \prod_{w \in D_\lambda} \left(\alpha \, (l(w) + 1) + \beta \, a(w)\right) \\
\intertext{and}
c'_\lambda(\alpha, \beta) & \coloneqq \prod_{w \in D_\lambda} \left(\alpha \, l(w) + \beta \, (a(w) + 1)\right).
\end{align*}
In particular, the product of the hook lengthes is given by
\[
h(\lambda) \coloneqq c_\lambda(1, 1) = c'_\lambda(1, 1) = \prod_{w \in D_\lambda} h(w).
\]

\paragraph{Characteristic classes.}

Let $\phi$ be a (rational) characteristic class (of complex vector bundles), i.e.~a polynomial in the Chern classes. We denote its homogeneous component of degree $k$ by $\phi_k$. In particular,
\[
\phi_k(F) \in H^{2 k}(X, \set Q)
\]
for any complex vector bundle $F$ on any manifold $X$.

\subsection{Equivariant cohomology}

\paragraph{Universal $G$-bundles.}
Let $G$ be a connected Lie group. Recall (e.g.~from~\cite{equivariant} and the references therein) that for each $n \in \set N_0$ there is an $n$-connected manifold $EG(n)$, on which $G$ acts freely. Let us fix this space throughout the discussion. The quotient is denoted by
\[
BG(n) \coloneqq EG(n)/G.
\]
The quotient map $EG(n) \to BG(n)$ is a principal $G$-bundle that is universal with respect to all manifolds of dimension at most $n$. Applying the universal property of $EG(n) \to BG(n)$ to the bundle $EG(m) \to BG(m)$ with $n \gg m$ induces differentiable maps $BG(m) \to BG(n)$ and $EG(m) \to EG(n)$ that are unique up to homotopy. Thus, the spaces $BG(n)$ and $EG(n)$, respectively, form a directed system. We set
\begin{align*}
EG & \coloneqq \colim_{n} EG(n) & BG \coloneqq \colim_{n} BG(n).
\end{align*}
The induced $G$-bundle $EG \to BG$ is a universal $G$-bundle for manifolds of arbitrary dimension.

\paragraph{Equivariant cohomology.}
Let $X$ be a $G$-manifold, i.e.~a manifold, on which $G$ acts. Then $G$ acts on $X \times EG$ freely by the diagonal action. The quotient is denoted by $X_G \coloneqq X \times_G EG$.
Then
\[
H^*_G(X, \set Q) \coloneqq H^*(X_G, \set Q)
\]
is the \emph{equivariant singular cohomology ring of $X$ (over the rationals $\set Q$)}. For example, 
$H^*_G(X, \set Q) = H^*(X, \set Q)$ if $G$ is the trivial group. For $n \gg m$, it is
\[
H^m_G(X, \set Q) = H^m(X \times_G EG(n), \set Q),
\]
i.e.~we can calculate each equivariant cohomology group of $X$ using $EG(n)$ instead of $EG$ with $n$ big enough.

\paragraph{Equivariant Borel--Moore homology.}

For a manifold $X$, one may define its Borel--Moore homology groups by
\[
H^{\BM}_*(X, \set Q) \coloneqq H_*(\hat X, \Set \infty, \set Q),
\]
where $\hat X$ is the one-point compactification of $X$. 

Let us return to the manifold $X$, on which $G$ acts. For each $m \in \set N_0$ one defines
\[
H^{G, \BM}_m(X, \set Q) \coloneqq H^\BM_{m + \dim EG(n)}(X \times_G EG(n), \set Q)
\]
for $n \gg m$. (Up to canonical isomorphism, this definition does not depend on $n$.)
This is called the \emph{$m$-th equivariant Borel--Moore homology group of $X$ (over the rationals)}. As we won't make use of ordinary singular homology (i.e.~homology with compact supports), we will drop the exponent ``$\BM$'' from now on.

As usual, the homology spaces $H^G_{-*}(X, \set Q)$ form naturally a module over the cohomology ring $H^*_G(X, \set Q)$ by the cap-product $\cap$.

\paragraph{Poincaré duality.}

Let $Y \subset X$ be a $G$-stable submanifold of dimension $m$ of $X$. By
\[
[Y]_G \in H^G_m(X, \set Q)
\]
we denote its equivariant fundamental class in $X$. This is defined by
\[
[Y]_G \coloneqq [Y \times_G EG(n)]
\]
for $n \gg 0$. In particular, there is an equivariant fundamental class $[X]_G \in H^G_{\dim X}(X, \set Q)$
of $X$ itself. Set
\[
D: H^*_G(X, \set Q) \to H^G_{\dim X - *}(X, \set Q),\ \alpha \mapsto \alpha \cap [X]_G.
\]
It follows from classical (i.e.~non-equivariant) Poincaré duality that this map is an isomorphism, called the \emph{equivariant Poincaré duality isomorphism}.

\paragraph{Proper push-forward.}

Let $f\colon X \to Y$ be a proper $G$-equivariant map between $G$-manifolds. For every $n \in \set N_0$, this induces a push-forward operator
\[
f_*\colon H_{* + \dim EG(n)}(X \times_G EG(n), \set Q) \to
H_{* + \dim EG(n)}(Y \times_G EG(n), \set Q).
\]
This in turn induces push-forward maps
\[
f_*\colon H^G_*(X, \set Q) \to H^G_*(Y, \set Q).
\]
In particular, one has $f_*([X]_G) = [f(X)]_G$ if $f$ is an embedding. By Poincaré duality,
we can define a map
\[
f_!\colon H^*_G(X, \set Q) \to H^{* + \dim Y - \dim X}(X, \set Q),\ \alpha \mapsto D^{-1} f_* D \alpha.
\]

\paragraph{Equivariant characteristic classes.}

Let $F$ be an equivariant vector bundle over $X$, i.e.~there a fibre-wise linear action of $G$ given on $F$ that is compatible with the action on $X$. Then $F_G \coloneqq F \times_G EG$ is canonically a vector bundle over $X_G$ of the same rank as $F$. Given a characteristic class $\phi$, one calls
\[
\phi_G(F) \coloneqq \phi(F_G) \in H^*_G(X, \set Q)
\]
the \emph{equivariant $\phi$-class} of $F$. Let us denote by $j\colon X \to X_G$ the inclusion of $X$ into $X_G$ as a fibre of $X_G \to BG$. Note that this induces a map
\[
j^*\colon H^*_G(X, \set Q) \to H^*(X, \set Q).
\]
It is $j^* F_G = F$, e.g.~by naturality of the characteristic classes, we have
\[
j^* \phi_G(F) = \phi(F).
\]
In particular, we can recover the non-equivariant characteristic classes of $F$ by its equivariant ones.

\paragraph{Localisation.}

Assume that the fixpoint locus $X^G$ of $X$ consists only of a finite number of isolated points. Let $i\colon X^G \to X$ be the inclusion map, which is proper. Thus it induces a $H^*(BG, \set Q)$-linear map
\[
i_!\colon H^*_G(X^G, \set Q) \to H^*_G(X, \set Q).
\]

Let us further assume that $G = T$ and that $X$ is a complex manifold, on which the given $T$-action can be extended to a holomorphic $\set C^\units$-action. We have $BT = \set P^\infty$. In particular,
\[
H^*(BT, \set Q) = \set Q[u]
\]
with $\deg u = 2$.
Given a $\set Q[u]$-module $M$, let us denote the localisation of $M$ with respect to the family $\Set{u^n}$ by $M'$. In this situation, the map $i_!$ induces an isomorphism
\[
i_!\colon H^*_T(X^T, \set Q)' \to H^*_T(X, \set Q)'
\]
of $\set Q[u, u^{-1}]$-modules.

This is the statement of the so-called localisation theorem (\cite{loc1}; and also~\cite{loc2}).

\paragraph{Calculation of equivariant characteristic classes.}

Let $\phi$ be a multiplicative characteristic class over the algebra $A$ (over $\set Q$), i.e.~there is a power series $f \in 1 + x A[[x]]$ such that
\[
\phi = \prod_{i} f(x_i)
\]
where the $x_i$ are the Chern roots. 

Again, we assume that $G = T$ and that $X$ is a complex manifold of complex dimension $n$, on which the given $T$-action can be extended to a holomorphic $\set C^\units$-action. Using the localisation formula, one can prove that
\begin{equation}
\label{localisation}
\phi(X) = \sum_{x \in X^T} \prod_{w \in W(x)} \frac{f(w u)} {w u} \, D^{-1} [x]_T,
\end{equation}
where $W(x)$ is the multiset of the $n$ weights of the induced $T$-action on the tangent space $T_x$.

\paragraph{Leray--Serre spectral sequence.}
Let us return to a general connected Lie group $G$.
As $G$ is connected, the exact homotopy sequence of the fibration $EG \to BG$ yields that $BG$ is simply-connected. It follows that the Leray--Serre spectral sequence associated to the locally trivial fibration $X_G \to BG$ with fibre $X$ yields a multiplicative spectral sequence with $E_2$-term
\[E_2^{*, *} = H^*(BG, \set Q) \otimes H^*(X, \set Q)\]
converging to $H^*_G(X, \set Q)$. Let us denote the filtration of $H^*_G(X, \set Q)$ given by the spectral sequence by $F^* H^*_G(X)$. The edge morphism $H^*_G(X, \set Q) \to H^*(X, \set Q)$ is simply the map $j^*$ induced by the inclusion $j\colon X \to X_G$ of $X$ into $X_G$ as a fibre of $X_G \to BG$.

\paragraph{A degenerate case.}
From now on assume that $X$ and $BG$ have no cohomology in odd degrees. For example, the latter is the case if $G = T$ as then $BG \simeq \set P^\infty$. Then the Leray--Serre spectral sequence degenerates already at the $E_2$-term, i.e.~there is a filtration $F^* H^*_G(X, \set Q)$ such that
\[
H^*(BG, \set Q) \otimes H^*(X, \set Q) = \gr^* H^*_G(X, \set Q)
\]
holds for the associated graded algebras. The filtration is given by
\[
F^p H^{p + q}_G(X, \set Q) = \Set{\alpha \cup \phi \mid \alpha \in H^p(BG, \set Q), \phi \in H^q_G(X, \set Q)}.
\]

In particular, the non-equivariant cohomology can be read off the equivariant one.

\subsection{The total space of $\sheaf O_{\set P^1}(-\gamma)$ on the projective line}

\paragraph{The cohomology.}

Let $\gamma \in \set Z$. Let $X(\gamma)$ be the total space of $\sheaf O_{\set P^1}(-\gamma)$. It is a quasi-projective complex surface homotopic to $\set P^1$. Let us denote by
\[
h = D^{-1}[L] \in H^2(X(\gamma), \set Q)
\]
the cohomology class of a fibre.

It follows that
\[
H^*(X(\gamma), \set Q) = H^*(\set P^1, \set Q) = \set Q[h]/(h^2).
\]

Let $Z$ be the zero section. Then
\[
[Z] = - \gamma \, [L],
\]
which follows, for example, from the following equations in the ordinary (compactly supported) homology:
\begin{align*}
D^{-1}[L] \cap [Z] & = [\pt] & \text{and} & & D^{-1}[Z] \cap [Z] & = - \gamma \, [\pt].
\end{align*}

\paragraph{Coordinates.}

We will regard $X(\gamma)$ as the quotient space
$((\set C^2 \setminus \Set 0) \times \set C)/\set C^\units$,
where the $\set C^\units$-action is given by
\begin{multline*}
\set C^\units \times ((\set C^2 \setminus \Set 0) \times \set C)
\to (\set C^2 \setminus \Set 0) \times \set C,\
\\
(s, (z_0, z_1, w)) \mapsto (s \, z_0, s \, z_1, s^{-\gamma} \, w).
\end{multline*}
The orbit of $(z_0, z_1, w)$ under this action shall be denoted by $[z_0:z_1:w^{- \frac 1 \gamma}]$.

\paragraph{The canonical class.}

The canonical sheaf of $X(\gamma)$ is given by
\[
\omega_{X(\gamma)} = p^* \left(\sheaf O(\gamma - 2)\right)
\]
where $p\colon X(\gamma) \to \set P^1$ is the projection map from the total space to the projective line. It follows that the canonical class $K_{X(\gamma)}$ of $X(\gamma)$ is given by
\[
K_{X(\gamma)} = (\gamma - 2) \, h.
\]
In particular, $K_{X(2)} = 0$.

In any case, we have that the Euler class $e_{X(\gamma)}$ of $X$ vanishes as $H^4(X(\gamma)) = 0$.

\paragraph{A $T$-action.}

There is a well-defined action
\[
T \times X(\gamma) \to X(\gamma),\ 
(t, [z_0:z_1:w^{-1/\gamma}]) \mapsto [t^{-1} \, z_0, t \, z_1, (t\, w)^{-1/\gamma}]
\]
of the Lie group $T$ on $X(\gamma)$ (see also~\cite{jack}). This action has two isolated fixpoints, namely
\begin{align*}
x_0 & \coloneqq [0:1:0] & \text{and} && x_1 & \coloneqq [1:0:0].
\end{align*}

Let $L_i$ be the fibre of $X(\gamma)$ that goes through $x_i$. Then $L_i$ is a $T$-invariant submanifold of $X$ with fixpoint $x_i$. The torus $T$ acts on $T_{x_0} L_0$ with weight $1$ and on $T_{x_1} L_1$ with weight $1 - \gamma$.

Furthermore, $Z$ is a $T$-invariant submanifold of $X(\gamma)$ with fixpoints $x_0$ and $x_1$. The weights of the induced $T$-action on $T_{x_0} Z$ and $T_{x_1} Z$ are given by $-1$ and $1$ respectively.

In particular, the weights of the $T$-action on $T_{x_0} X(\gamma)$ are given by $\Set{1, -1}$ and on $T_{x_1} X(\gamma)$ by $\Set{1 - \gamma, 1}$.

For all these statements, we refer to~(\cite{jack}).

\paragraph{Equivariant cohomology.}

We continue to denote the generator of $H^*(BT, \set Q)$ by $u$. As $X(\gamma)$ is a complex manifold with vanishing cohomology in odd degrees, we have
\[
H^*(X(\gamma), \set Q)[u] \simeq \gr^* H^*_G(X(\gamma), \set Q).
\]
Let us denote by $j\colon X(\gamma) \to X(\gamma)_T$ the inclusion of $X(\gamma)$ as a fibre. Then a basis of $H^*_G(X(\gamma), \set Q)'$ as a $\set Q[u, u^{-1}]$-module is given by $1$ and $D^{-1} [Z]_T$. It is
$j^* 1 = 1$ and $j^* D^{-1} [Z]_T = D^{-1} [Z]$.

Furthermore, we have by the localisation formula
\begin{align*}
%\eqref{l0}
[L_0]_T & = \gamma^{-1}(-[Z]_T + (1 - \gamma) \, u) \\
\intertext{and}
%\eqref{l1}
[L_1]_T & = \gamma^{-1}(-[Z]_T + u).
\end{align*}

\subsection{The Hilbert scheme of points on a surface}

\paragraph{The setting.}

Let $X$ be a quasi-projective complex surface. For each $n \in \set N_0$ we denote by $X^{[n]}$ its \emph{Hilbert scheme of zero-dimensional subschemes of length $n$}, i.e.~a morphism $T \to X^{[n]}$ from a quasi-projective complex manifold $T$ to $X^{[n]}$ is the same as giving a closed subscheme of $T \times X$ that is flat and finite of degree $n$ over $T$. As the Hilbert scheme is a fine moduli space, it comes with a universal subscheme $\Xi^n$ of $X^{[n]} \times X$, which is flat and finite of degree $n$ over $X^{[n]}$. Given a closed point $\xi$ of $X^{[n]}$, i.e.~a zero-dimensional subscheme of $X$ of length $n$, the fibre of $\Xi^n$ over $\xi$ is simply $\xi$ itself.

By a result of Fogarty, the quasi-projective scheme $X^{[n]}$ is smooth of dimension $2 n$.

\paragraph{An incidence variety.}

Given $\xi \in X^{[n]}$, we denote its support by
\[
\supp \xi \coloneqq \supp \sheaf O_\xi.
\]
This is a finite subset of $X$, whose cardinality lies between $1$ and $n$. If $\xi' \in X^{[m + n]}$ is another zero-dimensional subscheme of $X$ such that $\sheaf O_{\xi'}$ is an extension of $\sheaf O_\xi$, we write $\xi \subset \xi'$. In particular, $\supp \xi \subset \supp \xi'$. In this situation, the kernel of the map $\sheaf O_{\xi'} \to \sheaf O_{\xi}$ is denoted by $\sheaf I_{\xi, \xi'}$.

Now we can define the following incidence variety: Let $m \in \set N_0$. Let $\Xi^{m, n}$ be the subvariety of $X^{[m + n]} \times X^{[n]} \times X$ that is the closure of the subset
\[
\Set{(\xi', \xi, x) \in X^{[m + n]} \times X^{[n]} \times X\mid \supp \sheaf I_{\xi, \xi'} \subset \Set x}.
\]

We denote the projections of $X^{[m + n]} \times X^{[n]} \times X$ onto its factors uniformly by $p$, $q$ and $r$ as given in the diagram:
\[
\begin{CD}
X^{[m + n]} @<p<< X^{[m + n]} \times X^{[n]} \times X @>r>> X \\
& & @VqVV \\
& & X^{[n]}.
\end{CD}
\]
Note that the restriction $p|_{\Xi^{m, n}}\colon \Xi^{m, n} \to X^{[m + n]}$ is proper.

\paragraph{Nakajima's creation operators.}

Let $\alpha \in H^*(X)$. With the notions of the previous paragraph, we define
\begin{multline*}
q_m(\alpha)\colon H^*(X^{[n]}, \set Q) \to H^{* + \deg \alpha + 2 \,m - 2}(X^{[m + n]}, \set Q),\ 
\\
\beta \mapsto
D^{-1} p_*((r^* \alpha \cup q^* \beta) \cap [\Xi^{m, n}]).
\end{multline*}
Note the cohomological degree of the operator $q_m$, which is $\deg \alpha + 2\, m - 2$ (when we assume that $\alpha$ is homogeneous). This operator is called a \emph{creation operator}.

Given a partition $\lambda = (\lambda_1, \lambda_2, \dots)$ we define another operator $q_\lambda(\alpha)$ as follows: Let $r \coloneqq l(\lambda)$ the length of $\lambda$. The diagonal map $\delta^r\colon X \to X^r$ induces a proper push-forward map $\delta^r_!\colon H^*(X) \to H^*(X^r) = H^*(X)^{\otimes r}$. Let us write
\[
\delta^r_! \alpha = \sum \alpha_{i_1} \otimes \cdots \otimes \alpha_{i_r}.
\]
Then we define
\begin{multline*}
q_\lambda(\alpha)\colon H^*(X^{[n]}, \set Q) \to H^{* + \deg \alpha + 2 \abs\lambda + 2 l(\lambda) - 4}(X^{[n + \abs\lambda]}, \set Q)\ 
\\
\beta \mapsto \sum q_{\lambda_1}(\alpha_{i_1}) \, \cdots \, q_{\lambda_n}(\alpha_n)(\beta).
\end{multline*}

Let us denote by $\vac$ the unit in $H^0(X^{[0]}, \set Q) = \set Q$, the \emph{vacuum}. By results of I.~Grojnowski and H.~Nakajima, all elements in the rings $H^*(X^{[n]}, \set Q)$ can be created by applying creation operators on the vacuum $\vac$. More precisely, for every class $\alpha \in H^*(X^{[n]}, \set Q)$ and a basis $(\alpha_i)_i$ of $H^*(X, \set Q)$, there is a unique polynomial $Q$ in the $q_m(\alpha_i)$ such that
$\alpha = Q \vac$.

\paragraph{Characteristic classes.}

Let $\phi$ be a multiplicative characteristic class over a $\set Q$-algebra $A$. Then there are unique sequences $(a^\lambda_\phi)_\lambda$, $(b^\lambda_\phi)_\lambda$, $(c^\lambda_\phi)_\lambda$ and $(d^\lambda_\phi)_\lambda$ in $A$ where $\lambda$ runs through all partitions such that for all quasi-projective complex surfaces $X$, we have
\[
\sum_{n \ge 0} \phi(X^{[n]}) = \exp\sum_{\lambda} \left(a^\lambda_\phi \, q_\lambda(1) + b^\lambda_\phi \, q_\lambda(K_X) + c^\lambda_\phi \, q_\lambda(e_X) + d^\lambda_\phi \, q_\lambda(K^2_X)\right) \vac,
\]
where $K_X$ and $e_X$ denote the canonical and Euler class, respectively, of the surface $X$.

A closed formula for all coefficients is not known.

For the Chern character $\ch$, there is an analoguous result: There are unique sequences $(a^\lambda)_\lambda$, $(b^\lambda)_\lambda$, $(c^\lambda)_\lambda$ and $(d^\lambda)_\lambda$ in $\set Q$ such that for all quasi-projective complex surfaces $X$, we have
\begin{multline*}
\sum_{n \ge 0} \ch(X^{[n]}) \\
= \sum_{\lambda} \left(a^\lambda \, q_\lambda(1) + b^\lambda \, q_\lambda(K_X) + c^\lambda \, q_\lambda(e_X) + d^\lambda \, q_\lambda(K_X^2)\right) \, e^{q_1(1)} \vac.
\end{multline*}

\subsection{Combinatorial formulas}

\paragraph{The Lagrange--Good formula.}

Let $A$ be a ring and let $s \in \set N$. Let us write $z = (z_1, \dots, z_s)$.
Let $f_1, \dots, f_s \in A[z]$ be formal power series with $f_i \in z_i \set Q[[z]]$ for all $i = 1, \ldots, s$. Furthermore we assume that the $\left.\frac{\partial f_i}{\partial z_i}\right|_{z_i = 0}$ are invertible in $A$. We can then expand any Laurent series $g \in A((z))$ in terms of the $f_i(z)$, i.e.~
\[
g = \sum_{k \in \set Z^s} c_k f^k
\]
with $f^{(k_1, \dots, k_s)} \coloneqq f_1^{k_1} \, \cdots \, f_s^{k_s}$.

Then the following holds for the coefficients $c_k$:
\[
c_k = \res_z \frac{g(z)\, \frac{\partial f}{\partial z}}{f^{(k_1 + 1, \ldots, k_s + 1)}},
\]
where the residue is taken with respect to each variable $z_i$.

This formula is called the ``Lagrange--Good formula''. For details and a proof see, e.g.~\cite{langrange}.

\section{The Hilbert schemes of points on the total space of $\sheaf O_{\set P^1}(-2)$}

\subsection{On the cohomology}

\paragraph{The non-vanishing cohomology groups.}
Let $n \in \set N_0$. As the cohomology of $X(\gamma)^{[n]}$ is generated by Nakajima's creation operators, it has to be concentrated in degrees $0, 2, \dots, 2 n$, i.e.~$X(\gamma)^{[n]}$ has no cohomology in odd degrees and no cohomology in degrees greater than the degree of the middle cohomology group. A $\set Q$-basis of the cohomology group $H^{2 n}(X(\gamma)^{[n]}, \set Q)$ is given by
\[
q_{\lambda_1}(h) \cdots q_{\lambda_r}(h) \, \vac
\]
where $\lambda = (\lambda_1, \dots, \lambda_r)$ runs through all partitions of $n$.

\paragraph{Characteristic classes for $\gamma = 2$.}

Assume that $\gamma = 2$ in this paragraph. In particular, we have $K_{X(2)} = e_{X(2)} = 0$.

Let $r \in \set N$. Let $\delta^r\colon X \to X^r$ the diagonal embedding. Due to degree reasons, we have $\delta^r_! 1 = 0$ for $r \ge 3$. For $r = 2$, we have
\[
\delta^2_! 1 = - \frac 1 2 \, (D^{-1}[Z] \boxtimes D^{-1}[Z]) = - 2 \, (h \boxtimes h).
\]
It follows that for any multiplicative class $\phi$, we have
\[
\sum_{n \ge 0} \phi(X^{[n]}) = \exp\left(\sum_{k \ge 1} a^{(k)}_\phi \, q_k(1) - 2 \,
\sum_{k \ge l \ge 1} a^{(k, l)}_\phi \, q_k (h) q_l(h)\right)\, \vac.
\]
Due to degree reasons, this implies
\[
\sum_{n \ge 0} \phi_{n - 1}(X^{[n]}) = \left(\sum_{k \ge 1} a^{(k)}_\phi \, q_k(1)\right) \, \vac
\]
and
\[
\sum_{n \ge 0} \phi_{n}(X^{[n]}) = \exp\left(- 2 \, \sum_{k \ge l \ge 1} a^{(k, l)}_\phi \, q_k(h)
q_l(h)\right) \, \vac.
\]
In particular, we can read off the $a^{(k, l)}_\phi$ (that are universal for every surface!) once we know the classes $\phi_{n}(X^{[n]})$.

An analoguous result holds for the Chern character, namely
\[
\sum_{n \ge 0} \ch_{n}(X^{[n]}) = - 2 \, \left(\sum_{k \ge l \ge 1} a^{(k, l)} \,
q_k(h) q_l(h)\right) \, \vac.
\]
Again we can read off the $a^{(k, l)}$ once we know the classes $\ch_{n}(X^{[n]})$.

\subsection{The equivariant cohomology of the Hilbert scheme of points on $X(\gamma)$}

\paragraph{A $T$-action on $X(\gamma)^{[n]}$.}

Recall the action of the one-dimensional torus $T$ on $X(\gamma)$ we have defined above. By the universal property of the Hilbert scheme $X(\gamma)^{[n]}$ this action induces a $T$-action on $X(\gamma)^{[n]}$. This is the action we are interested in in the sequel. Consider a fixpoint $\xi \in X(\gamma)^{[n]}$. As $X(\gamma)$ has only two isolated fixpoints $x_0$ and $x_1$, it follows that $\supp \xi \in \Set{x_0, x_1}$. Moreover it turns out that $X(\gamma)^{[n]}$ has only isolated fixpoints and one can enumerate them by pairs $(\lambda^0, \lambda^1)$ of partitions such that $l(\lambda^0) + l(\lambda^1) = n$. This follows from the work of G.~Ellingsrud and S.~A.~Stømme,~\cite{fixpoints}; but see also~\cite{jack}. Let us denote the subscheme corresponding to the pair $(\lambda^0, \lambda^1)$ by $\xi_{\lambda^0, \lambda^1}$. The weights of the induced action of $T$ at the tangent space $T_{\xi_{\lambda^0, \lambda^1}}$ are given by $W_{\lambda^0}(-1, -1)$ and $W_{\lambda^1}(\gamma - 1, 1)$ by~\cite{jack}.

\paragraph{Generators of the equivariant cohomology ring.}
Note that by the localisation theorem, the localised equivariant cohomology ring of $X(\gamma)^{[n]}$ has the classes $[\xi_{\lambda^0, \lambda^1}]_T$ as generators over $\set Q[u, u^{-1}]$. We follow~\cite{jack} and set
\[
[\lambda^0, \lambda^1]_T \coloneqq \frac{(-1)^n}{c_{\lambda^0}(-1, -1) \, c_{\lambda^1}(\gamma - 1, 1)}
\, u^{-1} [\xi_{\lambda^0, \lambda^1}]_T \in H^{2n}_T(X(\gamma)^{[n]}, \set Q)'.
\]

\paragraph{Jack symmetric functions.}

Let $\Lambda$ be the ring of symmetric functions (in the variables $x_i$). By $\Lambda_n$ we denote its subspace consisting of the homogeneous symmetric functions of degree $n$.

Given a partition $\lambda$ and a parameter $\alpha$, we denote the Jack symmetric function for the partition $\lambda$ and the parameter $\alpha$ by $P^{(\alpha)}_\lambda$. (For the definition, we refer to I.~G.~Macdonald's book~\cite{symmetric}.) It is $P^{(\alpha)}_\lambda \in \Lambda_{\abs\lambda}$. The Jack symmetric functions to the parameter $\alpha = 1$ are simply the Schur symmetric functions, i.e.~
\[
P^{(1)}_\lambda = s_\lambda.
\]

Later, we will also need the power symmetric function
\[
p_n \coloneqq \sum_i x_i^n.
\]
Following~\cite{symmetric}, we set
\[
p_\lambda \coloneqq p_{\lambda_1} p_{\lambda_2} \cdots
\]
for a partition $\lambda = (\lambda_1, \lambda_2, \ldots)$. (This convention differs from the one used in~\cite{jack}.)

By~\cite{jack}, there is a canonical isomorphism
\[
\chi\colon H^{2n}_T(X(\gamma)^{[n]}, \set Q) \to (\Lambda \otimes \Lambda)_n
\]
mapping
$[\lambda^0, \lambda^1]_T$ to $P_{\lambda^0}^{(1)} \otimes P_{\lambda^1}^{(\gamma - 1)^{-1}}$.
In particular, for $\gamma = 2$ we have
\[
\chi([\lambda^0, \lambda^1]_T) = s_{\lambda^0} \otimes s_{\lambda^1}.
\]

\paragraph{From the equivariant cohomology to the non-equivariant one.}

Recall that we have the edge morphism $j^*\colon H^{2n}_T(X(\gamma)^{[n]}, \set Q) \to H^{2n}(X^{[n]}, \set Q)$.
This can be composed with the inverse of $\chi$, so that we get a morphism
\[
\psi \coloneqq j^* \circ \chi^{-1}\colon (\Lambda \otimes \Lambda)_n \to H^{2n}(X^{[n]}).
\]
By the discussion in~\cite{jack} following their theorem~3.4, this morphism can be explicitely described by
\[
\psi\colon p_{\lambda^0} \otimes p_{\lambda^1} \mapsto
q_{\lambda^0_1}(h) q_{\lambda^0_2}(h) \cdots \,
q_{\lambda^1_1}(h) q_{\lambda^1_2}(h) \cdots \, \vac.
\]

\subsection{Equivariant multiplicative characteristic classes of the Hilbert scheme of points on $X(\gamma)$}

\paragraph{The general formula.}

For the rest of this article, let $\phi$ be a multiplicative characteristic class over a $\set Q$-algebra $A$, i.e.~there is a power series $f \in 1 + x A[[x]]$ such that
\[
\phi = \prod_i f(x_i),
\]
where the $x_i$ are the Chern roots. 

By the localisation formula~\eqref{localisation} and the descriptions of the weights, we have the following equality in $H^{[2n]}_T(X(\gamma)^{[n]}, \set Q)$:
\[
\begin{split}
& \phantom{={}} \phi^T_{n}(X(\gamma)^{[n]})
\\ 
& = \sum_{\substack{\lambda^0, \lambda^1 \\ \abs{\lambda^0} + \abs{\lambda^1} = n}}
\begin{aligned}[t]
& \phantom{\cdot{}} (-1)^n \frac{[u^n] \left(\prod_{w \in W_{\lambda^0} \uplus W_{\lambda^1}} f(w u)\right)}
{c_{\lambda^0}(-1, -1) \, c'_{\lambda^0}(-1, -1) \, c_{\lambda^1}(\gamma - 1, 1) \,
c'_{\lambda^1}(\gamma - 1, 1)}
\\
& \cdot u^{-n} \, [\xi_{\lambda^0, \lambda^1}]_T
\end{aligned}
\\
& = \sum_{\lambda^0, \lambda^1}
\frac{[u^n] \left(\prod_{w \in W_{\lambda^0} \uplus W_{\lambda^1}} f(w u)\right)}
{c'_{\lambda^0}(-1, -1) \, c'_{\lambda^1}(\gamma - 1, 1)} \, [\lambda^0, \lambda^1].
\end{split}
\]
(For other formulas of this kind, see~\cite{boissiere} and~\cite{affine}.)

\paragraph{The specialisation to $\gamma = 2$.}

Let us now specialise to the case of a trivial canonical divisor, i.e.~$\gamma = 2$. Set $X := X(2)$. Then the previous formula can be simplified to
\[
\phi^T_{n}(X^{[n]}) = \sum_{\lambda^0, \lambda^1}
(-1)^{\abs{\lambda^0}} \frac{[u^n] \left(\prod_{w \in D_{\lambda^0} \uplus D_{\lambda^1}} f(h(w) u) \,
f(- h(w) u)\right)}
{h(\lambda^0) \, h(\lambda^1)}\, [\lambda^0, \lambda^1],
\]
where again the summation is over all pairs $(\lambda^0, \lambda^1)$ of partitions such that $\abs{\lambda^0} + \abs{\lambda^1} = n$.

\paragraph{The non-equivariant case.}

Recall that $\phi(X^{[n]}) = j^* \phi^T(X^{[n]})$. Thus we have
\begin{multline*}
\phi_{n}(X^{[n]}) = \sum_{\lambda^0, \lambda^1}
(-1)^{\abs{\lambda^0}} \frac{[u^n] \left(\prod_{w \in D_{\lambda^0} \uplus D_{\lambda^1}} f(h(w) u) \, f(- h(w) u)\right)}
{h(\lambda^0) \, h(\lambda^1)}\, \psi(s_{\lambda^0} \otimes s_{\lambda^1}),
\end{multline*}
as $\chi^{-1}([\lambda^0, \lambda^1]) = s_{\lambda^0} \otimes s_{\lambda^1}$.

\section{The proofs}

\subsection{Preparations}

\paragraph{Specialising the symmetric functions.}

We continue to use the definition $X = X(2)$ until the end of this article.

Let us define a map
\[
\rho\colon H^{2n}(X^{[n]}, \set Q) \to A[x, y],
\]
which is given by
\[
\rho\colon q_{\lambda_1}(h) \, \cdots \, q_{\lambda_n}(h) \vac \mapsto
(x^{\lambda_1} + y^{\lambda_1}) \, \cdots \, (x^{\lambda_n} + y^{\lambda_n}).
\]
The image of this map lies in the subspace of the polynomials of degree $n$. Recall the map
\[
\psi\colon (\Lambda \otimes \Lambda)_n \to H^{2n}(X^{[n]}, \set Q).
\]
We have
\[
\rho \circ \psi\colon (\Lambda \otimes \Lambda)_n \to A[x, y],\
(f_0 \otimes f_1) \mapsto f_0(x, y) \, f_1(x, y).
\]
Here, for a symmetric function $f \in \Lambda$, the expression $f(x, y)$ means to substitute $x$ for $x_1$, $y$ for $x_2$ and $0$ for $x_i$ with $i \ge 3$.

With the results of the previous subsection, we thus have
\begin{multline*}
\rho(\phi_{2n}(X^{[n]}))
\\
= \sum_{\lambda^0, \lambda^1}
(-1)^{\abs{\lambda^0}} \frac{[u^n] \left(\prod_{w \in D_{\lambda^0} \uplus D_{\lambda^1}} f(h(w) u) \, f(- h(w) u)\right)}
{h(\lambda^0) \, h(\lambda^1)}\, s_{\lambda^0}(x, y) \, s_{\lambda^1}(x, y).
\end{multline*}

\paragraph{First simplifications.}

Set
\[
F(x) \coloneqq f(x) \, f(-x) \in 1 + x \, A[[x]].
\]
Note that this is an even power series.

Our next task is to simplify the power series
\begin{multline*}
Z(x, y) \coloneqq \sum_{n \ge 0} \rho(\phi_{2n}(X^{[n]})) \\
= \sum_{\lambda^0, \lambda^1}
(-1)^{\abs{\lambda^0}} \frac{[u^n] \left(\prod_{w \in D_{\lambda^0} \uplus D_{\lambda^1}} f(h(w) u) \, f(- h(w) u)\right)}
{h(\lambda^0) \, h(\lambda^1)}\, s_{\lambda^0}(x, y) \, s_{\lambda^1}(x, y)
\\
= \sum_{\lambda^0, \lambda^1}
(-1)^{\abs{\lambda^0}} \frac{[u^n] \left(\prod_{w \in D_{\lambda^0} \uplus D_{\lambda^1}} F(h(w) u)\right)}
{h(\lambda^0) \, h(\lambda^1)}\, s_{\lambda^0}(x, y) \, s_{\lambda^1}(x, y)
\end{multline*}
in $A[[x, y]]$.

In order to do this, note that by the definition of the Schur polynomials (see, e.g.~\cite{symmetric}),
we have
\[
s_{\lambda}(x, y) = 0
\]
for $l(\lambda) \ge 3$ and
\[
s_{(a, b)}(x, y) = \frac{x^{a + 1} \, y^b - y^{a + 1} \, x^b}{x - y}
\]
for $a \ge b \ge 0$.

Therefore we can simplify the sum on the right hand side of the defining equation of $Z$ as all partitions of length three or greater do not contribute. The multiset of hook lengthes of a partition $\lambda$ of the form $(a, b)$, $a \ge b$ is given by
\[
\Set{1, \ldots, b; 1, \ldots, \widehat{a - b + 1}, \ldots, a + 1}.
\]
(The ``hat'' means to leave out that particular entry.)
In particular, we have
\[
h((a, b)) = \frac{(a + 1)!\, b!}{a - b + 1}.
\]

Thus we have
\begin{multline*}
Z = \frac 1{(x - y)^2} \, \sum_{\substack{a \ge b \ge 0 \\ c \ge d \ge 0}}^\infty
(-1)^{a + b} \frac{(a - b + 1)\, (c - d + 1)}{(a + 1)! \, b! \, (c + 1)!\, d!} \\
\cdot (x^{a + 1} y^b - y^{a + 1} x^b) \, (x^{c + 1} y^d - y^{c + 1} x^d)\\
\cdot [u^{a + b + c + d}]\frac{\prod_{w = - b}^{a + 1} F(w u) \, \prod_{w = - d}^{c + 1} F(w u)}{F((a - b + 1) u) F((c - d + 1) u)}.
\end{multline*}
By reindexing ($a + 1$ becomes $a$, $c + 1$ becomes $c$ and symmetrising between $a$ and $b$ and $c$ and $d$, respectively), the expression for $Z$ can also be written as
\begin{multline*}
Z = - \frac 1 {(x - y)^2} \, \sum_{a, b, c, d \ge 0}
(-1)^{a + b} \frac{(a - b) \, (c - d)}{a! \, b! \, c! \, d!} \, x^{a + c} \, y^{b + d}
\\
\cdot [u^{a + b + c + d - 2}]\frac{\prod_{w = - b}^{a} F(w u) \, \prod_{w = - d}^{c} F(w u)}{F((a - b) u) F((c - d) u)}.
\end{multline*}

Let us now collect all the coefficients in front of the monomials $x^r y^s$. This yields:
\begin{multline*}
Z = - \frac 1 {(x - y)^2} \, \sum_{r, s \ge 0} x^r y^s \,
[u^{r + s - 2}]
\sum_{a = 0}^r \sum_{b = 0}^s (-1)^{a + b}
\frac{(a - b) \, (r + s - (a + b))}{a! \, b! \, (r - a)! \, (s - b)!} \\
\cdot
\frac{\prod_{w = a - r}^{a} F(w u) \, \prod_{w = b - s}^{b} F(w u)}{F((a - b) u) F((r + s - (a + b)) u)}.
\end{multline*}

\paragraph{Some combinatorics.}

The calculations in this paragraph mimick those in~\cite{affine}.
Let $k \in \set N_0$. Recall the beginning of the classical summation formula
\[
\sum_{w = 1}^s w^k = \frac{s^{k + 1}}{k + 1} + \frac{s^k}{2} + \ldots,
\]
where ``$\ldots$'' denotes termes of lower degree in $s$.

Let $F \in 1 + x A[[x]]$ be a power series with constant term one. We write
\[
\log F = \sum_{k \ge 1} L_k x^k.
\]
Let $r, s, a, b \in \set N_0$ and $a \leq r$ and $b \leq s$. By the summation formula, we have
\begin{multline*}
\frac{\prod_{w = a - r}^{a} F(w u) \, \prod_{w = b - s}^{b} F(w u)}{F((a - b)  u) F((r + s - (a + b))  u)} \\
= \exp\left(\sum_{k = 1}^\infty L_k \, u^k \, \left(\sum_{w = 1}^a w^k + \sum_{w = 1}^{r - a} (-w)^k
+ \sum_{w = 1}^b w^k + \sum_{w = 1}^{s - b} (-w)^k
\right.\right.
\\
\left.\left. 
- (a - b)^k - (r + s - (a - b))^k\right)\right) \\
= \exp\left(\sum_{k = 1}^\infty L_k u^k \, \left(
(r + 1) \, a^k + (s + 1) \, b^k - (a - b)^k - (b - a)^k + \ldots
\right)\right)
\\
= \frac{F^{r + 1}(a u) \, F^{s + 1}(b u)}{F((a - b) u) \, F((b - a) u} + \ldots,
\end{multline*}
where ``$\ldots$'' collects all terms in which the total degree in $a$ and $b$ is less than the total degree in $u$.

In what follows, we continue to use that notion of the symbol ``$\ldots$''.

We use the following lemma from~\cite{affine}:
\[
\sum_{a = 0}^r \frac{(-1)^a a^k}{a! \, (r - a)!} = \begin{cases}
0 & \text{for $k < r$}
\\
(-1)^r & \text{for $k = r$}
\end{cases}
\]
for any $r, k \in \set N_0$ and $k \leq r$.

Thus, we have
\begin{multline*}
[u^{r + s - 2}] \sum_{a = 0}^r \sum_{b = 0}^s (-1)^{a + b} \frac{(a - b)((r - s) - (a - b))}
{a! \, b! \, (r - a)! \, (s - b)!} \\
\cdot \left(\frac{F^{r + 1}(a u) \, F^{s + 1}(b u)}{F((a - b) u) \, F((b - a) u)} + \ldots\right)
\\
= [u^{r + s}] \sum_{a = 0}^r \sum_{b = 0}^s (-1)^{a + b} \frac{(a - b) (b - a)  u^2}
{a! \, b! \, (r - a)! \, (s - b)!} \, \frac{F^{r + 1}(a u) \, F^{s + 1}(b u)}{F((a - b) u) \, F((b - a) u)} + \ldots
\\
= (-1)^{r + s} [a^r b^s] \frac{F^{r + 1}(a) \, F^{s + 1}(b)}{F(a - b) \, F(b - a)} \, (a - b) \, (b - a).
\end{multline*}

\paragraph{Returning to $Z$.}

With the results of the previous paragraph, we now have the following expression for $Z$:
\[
Z(x, y) = - \frac 1 {(x - y)^2} \sum_{r, s = 0}^\infty (-x)^r \, (-y)^s \, [a^r b^s]
\frac{F^{r + 1}(a) \, F^{s + 1}(b)}{F(a - b) \, F(b - a)} \, (a - b) \, (b - a).
\]
Set
\[
G(z) \coloneqq \frac z {F(z)}.
\]
This is an odd power series, which is invertible with respect to the composition.
It is
\[
Z(x, y) = - \frac 1 {(x - y)^2} \sum_{r, s = 0}^\infty (-x)^r \, (-y)^s \, \res_{(a, b)}
\frac{G(a - b) \, G(b- a)}{G^{r + 1}(a) \, G^{s + 1}(b)}.
\]
By the Lagrange--Good formula, we have
\[
Z(-G(x), - G(y)) = \frac 1 {G'(x) \, G'(y)} \,
\frac{G(x - y)}{G(x) - G(y)} \, \frac{G(y - x)}{G(y) - G(x)}.
\]
As $G$ is an odd power series, this gives
\[
Z(G(x), G(y)) = \frac 1 {G'(x) \, G'(y)} \, \left(
\frac{G(x - y)}{G(x) - G(y)}\right)^2.
\]
We can reformulate this as follows: Let $g \in xA[[x]]$ be the compositional inverse of $G$. Then
\begin{equation}
\label{zeq}
Z(x, y) =  g'(x) g'(y) \left(
\frac{G(g(x) - g(y))}{x - y}
\right)^2.
\end{equation}

\subsection{The proof of the theorem and the corollary}

\paragraph{The proof of the theorem.}

Recall from the beginning that we have to prove the following:
\begin{equation}
\label{genshape}
\sum_{n \ge 0} \phi(X^{[n]}) = \exp\left(\sum_{k = 1}^\infty a_k \, q_k(1) +
\sum_{\substack{k, l = 1}}^\infty a_{k, l} \, q_k(h) q_l(h)\right) \, \vac,
\end{equation}
where the sequences $(a_k)_k$ and $(a_{k, l})_{k, l}$ in $A$ are defined by
\[
\sum_{k = 1}^\infty a_k \, k x^k = g(x)
\]
and
\[
\sum_{k, l \ge 1} a_{k, l} \, x^k y^l = 
\log\frac{g(x) - g(y)}{(x - y) \, (f(g(x) - g(y)) \, f(g(y) - f(g(x))))}.
\]

As we have already seen, there exist unique $a_k$ and $a_{k, l}$ with $a_{k, l} = a_{l, k}$, for which the formula~\eqref{genshape} holds true. It remains to calculate them.

As stated above, due to degree reasons, it is
\[
\sum_{n \ge 0} \phi_{n - 1}(X^{[n]}) = \left(\sum_{k \ge 1}^\infty a_k q_k(1) \, \right) \vac.
\]
These coefficients, however, were already calculated in~\cite{affine}. The main result of that article is
\[
\sum_{k = 1}^\infty a_k \, k x^k = g(x).
\]

Analoguously, we have
\[
\sum_{n \ge 0} \phi_{n}(X^{[n]}) = \exp\left(\sum_{k, l \ge 1} a_{k, l} \, q_k(h) q_l(h)\right) \, \vac.
\]

In order to finish the proof of the theorem, it remains to calculate the $a_{k, l}$: We have
\[
\exp\left(\sum_{k, l \ge 1} a_{k, l} \, (x^k + y^k) (x^l + y^l)\right)
= \rho\left(\sum_{n \ge 0} \phi_{2 n}(X^{[n]})\right)
= Z(x, y),
\]
i.e.~
\[
\sum_{k, l \ge 1} a_{k, l} \, (x^k + y^k) (x^l + y^l)
= \log Z(x, y).
\]
This yields
\begin{multline*}
\sum_{k, l \ge 1} a_{k, l} x^k \, y^l = \frac 1 2 \left(\log Z(x, y) - \log Z(x, 0) - \log Z(0, y)\right)
\\
= \log\frac{g(x) - g(y)}{(x - y) \, (f(g(x) - g(y)) \, f(g(y) - f(g(x))))},
\end{multline*}
which proves the theorem.

\paragraph{The proof of the corollary.}

We have to prove the following:
\[
\sum_{n \ge 0} \ch(X^{[n]}) = \left(\sum_{k \ge 1} a_k q_k(1) + \sum_{k, l \ge 1} a_{k, l} \, q_k(h) q_l(h)\right) \,
e^{q_1(1)} \, \vac,
\]
where the sequences $(a_k)_k$ and $(a_{k, l})_{k, l}$ in $\set Q$ are defined by
\[
\sum_{k \ge 1}^\infty a_k \, x^k = \sum_{m \ge 0} \frac{2}{(2 m + 1)!} x^{2 m + 1}
\]
and
\[
\sum_{k, l \ge 1} a_{k, l} \, x^k \, y^l
= \sum_{m \ge 0} \frac 2 {(2 m)!} \sum_{k + l = 2 \, m}
\left(1 - (-1)^k \binom{2 \, m}{k}\right) \, x^k \, x^l.
\]

From the general results in~\cite{boissiere} and~\cite{affine} it is again clear that a formula of the form given above exists indeed. It remains to calculate the coefficients. For the $a_k$, this has already been done in~\cite{boissiere}, i.e.~we have in fact
\[
\sum_{k \ge 1}^\infty a_k \, x^k = \sum_{m \ge 0} \frac{2}{(2 m + 1)!} x^{2 m + 1}.
\]

It remains to prove the expression given for the $a_{k, l}$. By degree reasons, it is
\[
\ch_{n} (X^{[n]}) = \left(\sum_{k + l = n} a_{k, l} q_k(h) q_l(h)\right) \, \vac
\]
for each $n \in \set N$. 

As $X^{[n]}$ is a (non-compact) symplectic manifold (as $X$ is), it is $\ch_{n} = 0$ for odd $n$. 
Thus, we may restrict to the case of even $n$.

Let $A := \set Q[\epsilon]/(\epsilon^2)$ be the ring of dual numbers of $\set Q$.
Fix $n = 2 m \in 2 \set N$ and consider the $A$-valued multiplicative characteristic class $\phi$ with
\[
\phi = \prod_i f(x_i),
\]
where
\[
f(x) = 1 + \epsilon \, x^n \in 1 + x \set A[[x]].
\]
Then
\[
n! \, \ch_{n} (X^{[n]}) = [\epsilon] \phi(X^{[n]}).
\]
Set
\[
G(x) \coloneqq \frac{x}{f(x) \, f(-x)} = x - 2 \epsilon \, x^{n + 1}.
\]
Then
\[
g(x) := x + 2 \epsilon \, x^{n + 1}
\]
is the multiplicative inverse of $G$.
A quick calculation yields
\begin{multline*}
[\epsilon]\log\frac{g(x) - g(y)}{(x - y) \, (f(g(x) - g(y)) \, f(g(y) - g(x)))}
\\
= 2 \left(\frac{x^{n + 1} - y^{n + 1}}{(x - y)} -  (x - y)^n\right).
\end{multline*}

So if
\[
\sum_{k, l \ge 1} a_{k, l} \, x^k \, y^l = 2 \, \sum_{m \ge 0} \frac{1}{(2m)!} \sum_{k + l = 2 \, m}
\left(1 - (-1)^k \binom {2n} k\right) \, x^k y^l,
\]
we have, by the theorem applied to $\phi$, that
\[
\ch_{n} (X^{[n]}) = \frac 1 {(2n)!} \, [\epsilon] \phi_{n}(X^{[n]}) = 
\left(\sum_{k + l = n} a_{k, l} \, q_k(h) q_l(h)\right) \vac,
\]
which proves the corollary.

\section{Addendum}

\subsection{Characteristic classes of tautological sheaves}

\paragraph{Tautological sheaves.}

Let $X$ be an arbitrary complex surface. Given a vector bundle $F$ on $X$, on defines
\[
F^{[n]} := p_*(r^* F \otimes \sheaf O_{\Xi^n}),
\]
where $p$ and $r$ are the projections from $X^{[n]} \times X$ onto its factors. As $\Xi^n$ is flat and finite of degree $n$ over $X^{[n]}$, it follows that $F^{[n]}$ is a vector bundle on $X^{[n]}$ of rank $n \cdot \rk F$. It is the \emph{tautological vector bundle on $X^{[n]}$ associated to $F$}.

\paragraph{Characteristic classes of $\sheaf O_X^{[n]}$ on $X(2)^{[n]}$.}
In particular, there is the vector bundle $\sheaf O_X^{[n]}$ of rank $n$ on $X^{[n]}$. As this is god-given as the tangent bundle on the Hilbert schemes, one may likewise ask for the value of characteristic classes on the $\sheaf O_X^{[n]}$ (see~\cite{affine} for example).

The methods used in this article also work in that case. Instead of giving full proofs, which are similar enough to the ones given here, we only state the result. Again, we specialise to the case $X = X(2)$, the total space of the line bundle $\sheaf O_{\set P^1}(-2)$.

\paragraph{Theorem.}

Let $\phi$ be multiplicative class of the $\set Q$-algebra $A$ defined by
\[
\phi = \prod_i f(x_i)
\]
for $f \in 1 + x A[[x]]$, where the $x_i$ are the Chern roots. Let $g \in xA[[x]]$ be the compositional inverse of
\[
\frac{x}{f(-x)}.
\]
Then, the multiplicative class evaluated on the tautological bundles $\sheaf O_X^{[n]}$ on the $X^{[n]}$ is given by
\[
\sum_{n \ge 0} \phi(\sheaf O_X^{[n]}) = \exp\left(\sum_{k \ge 1} a_k\, q_k(1) + \sum_{k, l \ge 1} a_{k, l} \,
q_k(h) q_l(h)\right)\, \vac,
\]
where the sequences $(a_k)_k$ and $(a_{k, l})_{k, l}$ in $A$ are defined by
\[
\sum_{k \ge 1} k a_k \, x^k = g(x)
\]
and
\[
\sum_{k, l \ge 1} a_{k, l} \, x^k y^l = \log\frac{x \, y \, (g(x) - g(y))}{(x - y) \, g(x) \, g(y)}.
\]

\subsection{Calculations}

\paragraph{The Chern character.}

Recall the definition of the $a^{(k, l)}$ involved in the universal formula of the Chern character of the tangent bundles of the Hilbert schemes of points on surfaces. Our formulas yield the following table:
\[
\begin{array}{|c||c||c|c||c|c|c|	}
\hline
(k, l) & (1, 1) & (3, 1) & (2, 2) & (5, 1) & (4, 2) & (3, 3)
\\
\hline
a^{(k, l)} & -\frac 3 2 & - \frac 5 {12} & \frac 5 {24} & - \frac 7 {360} & \frac 7 {180} & - \frac {7}{240}
\\
\hline
\end{array}
\]
(Recall that $a^{(k, l)} = 0$ for $k + l$ odd.)

\paragraph{The total Chern class.}

Let $c$ denote the total Chern class. Recall the definition of the $a^{(k, l)}_c$ involved in the universal formula for the total Chern class of the tangent bundles of the Hilbert schemes of points on surfaces. Our formulas yield the following table:
\[
\begin{array}{|c||c||c|c||c|c|c|	}
\hline
(k, l) & (1, 1) & (3, 1) & (2, 2) & (5, 1) & (4, 2) & (3, 3)
\\
\hline
a^{(k, l)} & \frac 3 2 & - 1 & - \frac 7 {4} & 2 & 2 & 3
\\
\hline
\end{array}
\]
(Again, it is $a^{(k, l)}_c = 0$ for $k + l$ odd.) 

\appendix


\begin{thebibliography}{99}

\bibitem{loc1}
Atiyah, M.; Bott, R.
\emph{The moment map and equivariant cohomology.}
Topology 23, 1--28 (1984).

\bibitem{boissiere}
Bossière, S.
\emph{Chern classes of the tangent bundle on the Hilbert scheme of points on the affine plane.}
J.~Algebr.~Geom.~14, No.~4, 761--787 (2005).

\bibitem{affine}
Boissière, S.; Nieper-Wißkirchen, M.
\emph{Universal formulas for characteristic classes on the Hilbert schemes of points on surfaces.}
arXiv:math.AG/0507470.

\bibitem{equivariant}
Brion, M.
\emph{Equivariant cohomology and equivariant intersection theory.}
arXiv:math.AG/9802063.

\bibitem{loc2}
Edidin, D.; Graham, W.
\emph{Localization in equivariant intersection theory and the Bott residue formula.}
Am.~J.~Math.~120, No.~3, 619--636 (1998).

\bibitem{fixpoints}
Ellingsrud, G.; Strømme, S.~A.
\emph{On the homology of the Hilbert scheme of points in the plane.}
Invent.~Math.~87, 343--352 (1987).

\bibitem{smoothness}
Fogarty, J.
\emph{Algebraic families on an algebraic surface.}
Am.~J.~Math.~90, 511--521 (1968).

\bibitem{multiplicative}
Hirzebruch, F. \emph{Topological methods in algebraic geometry.}
Classics in Mathematics. Berlin: Springer-Verlag. ix, 234 p.~(1995).

\bibitem{hilbg}
Grojnowksi, I. \emph{Instantons and affine algebras. I: The Hilbert scheme and vertex operators.}
Math.~Res.~Lett.~3, No.~2, 275--291 (1996).

\bibitem{langrange}
Krattenthaler, Ch. \emph{Operator methods and Lagrange inversion: A unified approach to Lagrange formulas.} Trans.~Am.~Math.~ Soc.~305, No.~2, 431--465 (1988).

\bibitem{jack}
Li, W.-P.; Qin, Zh.; Wang, W.
\emph{The cohomology rings of Hilbert schemes via Jack polynomials.}
Hurtubise, Jacques (ed.) et al., Algebraic structures and moduli spaces. Proceedings of the CRM workshop, Montréal, Canada, July 14--20, 2003. Providence, RI: American Mathematical Society (AMS). CRM Proceedings \& Lecture Notes 38, 249--258 (2004).

\bibitem{symmetric}
Macdonald, I.~G.
\emph{Symmetric functions and Hall polynomials. 2nd ed.}
Oxford Science Publications. Oxford: Clarendon Press.~x, 475 p.~(1998).

\bibitem{hilbn}
Nakajima, H. \emph{Heisenberg algebra and Hilbert schemes of points on projective surfaces.}
Ann.~Math.~(2) 145, No.~2, 379--388 (1997).

\bibitem{jackn}
Nakajima, H. \emph{Jack polynomials and Hilbert schemes of points on surfaces}
arXiv:math.AG/9610021.

\bibitem{vasserot}
Vasserot, E. \emph{Sur l'anneau de cohomologie du schéma de Hilbert de $\set C^2$.}
C.~R.~Acad.~Sci., Paris, Sér.~I, Math.~332, No.~1, 7--12 (2001).

\end{thebibliography}
\end{document}